\documentclass[11pt]{amsart}
\usepackage{amsmath,amscd,amssymb,amssymb}

\theoremstyle{plain}
\newtheorem{theorem}{Theorem}[section]
\newtheorem{prop}[theorem]{Proposition}

\newtheorem{cor}[theorem]{Corollary}

\theoremstyle{definition}

\newcommand{\N}{{\mathbb N}}

\newcommand{\HH}{{\mathcal H}}
\newcommand{\FF}{{\mathcal F}}

\newcommand{\kk}{\mathbf k}
\newcommand{\G}{{G}}

\newcommand{\GL}{{GL}}
\newcommand{\Stab}{\mathrm{Stab}}
\newcommand{\Ind}{\mathrm{Ind}}
\newcommand{\Hom}{\mathrm{Hom}}

\newcommand{\PP}{{P}}

\newcommand{\End}{\text{End}}

\newcommand{\F}{\mathtt{F}}

\newcommand{\M}{M_2}
\newcommand{\val}{v}

\newcommand{\la}{\lambda}

\newcommand{\C}{\mathbb{C}}
\newcommand{\inv}{^{-1}}
\newcommand{\bsl}{\backslash}


\begin{document}

\title{A note on Bruhat decomposition of $GL(n)$ over\\ local principal ideal rings}
\author{Uri Onn, Amritanshu
Prasad and Leonid Vaserstein} \maketitle \markboth{\textsc{URI ONN,
AMRITANSHU PRASAD AND LEONID VASERSTEIN}}{\textsc{BRUHAT
DECOMPOSITION OVER LOCAL PRINCIPAL IDEAL RINGS}}
\begin{abstract}
Let $A$ be a local commutative principal ideal ring. We study the double coset space of $\GL_n(A)$ with respect to the subgroup of upper triangular matrices.
Geometrically, these cosets describe the relative position of two
full flags of free primitive submodules of $A^n$. We introduce some invariants of the double cosets. If $k$ is the
length of the ring, we determine for which of the pairs $(n,k)$ the
double coset space depends on the ring in question. For $n=3$, we
give a complete parametrisation of the double coset space and
provide estimates on the rate of growth of the number of double
cosets.
\end{abstract}

\section{Introduction}

Let $A$ be a local principal ideal commutative ring and let
$\wp=(\pi)$ denote its maximal ideal. Denote by $k$ the length of
the ring, that is, the least $k$ such that $\wp^k=0$ ($k$ might be
infinite). Let $B$ denote the subgroup of upper triangular
matrices in $\G=\GL_n(A)$, the group  of invertible matrices with
entries in $A$. This paper concerns the description of the
double coset space $B \bsl \G / B$. Since $B$ is the stabiliser of the
standard flag in $A^n$, this space corresponds to the possible
relative positions of two flags that are isomorphic to the standard
flag (these are the flags whose reductions modulo $\wp$ are full flags
of vector spaces in $(A/\wp)^n$). We refer to such a flag as a \emph{full
free primitive flag over $A$ in $A^n$}. If
\begin{align*}
\FF&=\text{the space of full free primitive flags over $A$ in $A^n$} \\
B&=\Stab_{\G}(\text{standard flag})
\end{align*}
then one has
\[
\FF \times_\G \FF \longleftrightarrow B \backslash \G /B~.
\]
A complete description of the latter space for all $n$ contains the
embedding problem of pairs of $A$-modules, which is known to be of
\emph{wild type} in general \cite{RS,Schmidmeier}. Thus, one does
not expect a reasonably closed solution, and we aim at the more
modest goal of constructing some invariants of the double coset
space, describing the relations between them and exploring to what
extent they can distinguish double cosets. When $k=1$, in which case
the ring $A$ is a field, Bruhat decomposition says that $B\bsl G/B$
is parameterised by the symmetric group. In particular, the number
of double cosets does not depend on the field. The first natural
question is: for which pairs $(n,k)$ is the parametrisation
independent of the ring? The answer to this question is given in
Theorem~\ref{theorem:cases}.
The first two instances in which a dependence on the ring occurs are
$(3,k)$ with $k>2$ and $(4,2)$. For $(3,k)$ we give a complete
parametrisation of the double coset space and give estimates for its
size when $A$ is a finite ring (Section~\ref{sec:n_is_3}).

\subsection{Related problems}

Let $P_1,P_2 < G$ be finite groups. Let
$\rho_i=\Ind_{P_i}^{G}1=\C[G/P_i]$ be the representation of
$G$ induced from the trivial representation of $P_i$ over $\C$. The module
of intertwining operators $\Hom_G(\rho_1,\rho_2)$ can be identified
with the subalgebra $\C[P_2\bsl G/P_1]$ of left-$(P_2,P_1)$-invariant elements
in the group algebra $\C[G]$ via the map $\C[P_2\bsl G/P_1]\to \Hom(\rho_1,\rho_2)$ given by
\begin{equation*}
  f\mapsto T_f \mbox{ for each } f\in \C[P_2\bsl G /P_1],
\end{equation*}
where
\begin{equation*}
  T_f h (x) = \sum_{g\in G} h(xg)f(g) \mbox{ for each } h\in \rho_1=\C[G/P_1].
\end{equation*}
Let $A_i=A/\wp^i$ for $i \in \N$ be the inverse system of the
finite length quotients. Isomorphism types of finitely generated
$A_i$-modules correspond to Young diagrams with height bounded by
$i$.
The Young diagram given by $\la=(\la_1,\ldots,\la_j)$, with
$i\geq \la_1\geq \cdots \geq \la_j \geq 0$ corresponds to the
$A$-module
\[ M_{\la}=\oplus_{r=1}^{j} A/\wp^{\la_r}.\]
Let $\G_i=\GL_n(A_i)$,
$P_{\la}$ be the stabiliser of a submodule of type $\la$ in $A_i^n$,
and $B_i$ be the stabiliser of a full  flag of free submodules in
$A_i^n$.

Assume now that the residue field of $A$ is finite of order $q$. The induced
representations $\rho_i=\Ind_{B_i}^{\G_i}1$ play a significant role
in the representation theory of the groups $\G_i$ \cite{HG1,HG4} in
analogy with the role played by $\rho_1$ in the representation
theory of $\GL_n$ over finite fields \cite{Zelevinsky}. The latter
representation is studied in terms of the Hecke algebra
$\HH_{A,1}$, where
\[
\HH_{A,i}=\End_{\G_i}(\rho_i) \simeq \C[B_i \bsl \G_i /B_i]
\]
and one has $\HH_{A,1} \simeq \C[S_n]$, independent of the
characteristic of the residue field. The algebra $\HH_{A,i}$
continues to play an important role for $i > 1$, however, its structure depends on the characteristic of
the residue field.
As a starting point, one would like to know its
dimension - hence parameterise the double coset space $B_i \bsl \G_i
/B_i$.
The number of double cosets depends on the residue field unless $n \leq 2$ or $i=1$ or $n=3$ and $i=2$ (see Theorem~\ref{theorem:cases}).

Broadening the frame a bit, it is natural to consider the category
of diagrams over finite length $A$-modules; For a quiver $Q$ let
$\mathcal{D}_{Q}=\text{Fun}(Q, A\text{-mod})$ be the category of
functors from $Q$ (considered as a category) to the category of
finitely generated $A$-modules. These categories occur naturally in the study
of the above representations, in particular as the underlying sets
for the modules of intertwining operators between representations.
For example \cite{BO1}
\begin{align*}
&\PP_{i^m} \bsl \G_i / \PP_{i^m}& &\longleftrightarrow
\{\text{isomorphism types
of submodules of $M_{i^m}=A_i^m$}\} \hookrightarrow \text{Isom}(\mathcal{D}_{\{\bullet\}}) &&\\
&\PP_{i^m} \bsl \G_i / \PP_{\la}& &\longleftrightarrow
\{\text{isomorphism types
of pairs $(N \subseteq N') \subseteq M_{\la}$}\} \hookrightarrow \text{Isom}(\mathcal{D}_{\{\bullet \rightarrow \bullet\}}) &&\\
&\PP_{i^m} \bsl \G_i / B_i& &\longleftrightarrow \{\text{isomorphism
types of chains submodules of $A_i^m$}\} \hookrightarrow
\text{Isom}(\mathcal{D}_{\{\bullet \rightarrow \bullet \rightarrow
\cdots \rightarrow \bullet\}}) &&
\end{align*}
The category $\mathcal{D}_{\{\bullet\}}$ is nothing but the category
of $A$-modules. The full subcategory of $\mathcal{D}_{\{\bullet
\rightarrow \bullet\}}$ which consists of embeddings is of wild type
\cite{RS,Schmidmeier}. The full subcategory of
$\mathcal{D}_{\{\bullet \rightarrow \bullet \rightarrow \cdots
\rightarrow \bullet\}}$ which consists of embeddings is discussed in
\cite{Simson}.

\subsection{Notations} Throughout this paper $\pi$ denotes a generator of the maximal
ideal $\wp$ in $A$ and the order of $\pi$ is denoted by $k$, and might
be infinite (as in the preceding paragraphs). The valuation of a non-zero element $x\in A$ is denoted by $\val(x)$. For convenience we also write $\val(0)=k$. For any ideal $I$ in $A$,
$\val(I)$ is the semigroup of valuation values of elements in $I$.
The residue field of $A/\wp$ is denoted by $\kk$, $A^{\times}$ denotes the
multiplicative group, and $A_i=A/\wp^i$.

\subsection{Acknowledgements} We thank the Tata Institute of Fundamental Research, in particular Ravi Rao and Dipendra Prasad, for bringing us together and for their warm hospitality, which resulted in this manuscript.
We thank Markus Schmidmeier for supplying detailed information on
the embedding problem. The first author also thanks Uri Bader for
discussions which partly motivated this work.

\section{Invariants}
\subsection{Upper triangular row and column operations}
\label{sec:upper_triangular} The determination of the double coset space $B\bsl G/B$
is part of the more general question of determining the orbits of the
action of invertible upper triangular matrices by left and right multiplications on the set $M_{nm}$ of $n\times m$ matrices with entries in $A$.
When two matrices $\alpha$ and $\alpha'$ lie in the same orbit, we write $\alpha\sim \alpha'$.
Let $R_i$ denote the $i$th row and $C_j$ denote the $j$th column of a matrix.
Then two matrices lie in the same orbit if one can be obtained from the other by the following types of row and column operations:
\begin{itemize}
\item Multiplication of rows or columns by scalars (action of the \lq torus\rq)
\[
 R_i\to aR_i, \mbox{ with } {a\in A^\times} \mbox{ and } C_j\to aC_j, \mbox{ with } {a\in A^\times}.
\]
\item Addition of certain rows/columns to others (action of the \lq unipotent subgroup\rq)

\begin{align*}
&R_i\to R_i+\sum_{i'>i}a_{i'}R_{i'} \mbox{ with } a_{i'}\in A
\\ &C_j\to C_j+\sum_{j'<j}a_{j'}C_{j'} \mbox{ with } a_{j'}\in A
\end{align*}
\noindent
Note that only row (respectively, column) operations with $i' > i$ (respectively, $j' < j$) are allowed.
\end{itemize}
Since row operations commute with column operations and scaling
operations normalise addition operations, two matrices in
$M_{nm}$ lie in the same orbit if and only if one can
be obtained from the other by a sequence of scaling operations
followed by row operations $R_i\to R_i+\sum_{i'>i}a_{i'}R_{i'}$ with
$i$ increasing from $1$ to $n$ followed by column operations $C_j\to
C_j+\sum_{j'<j}a_{j'}C_{j'}$ with $j$ decreasing from $n$ to $1$.

\subsection{Decomposability} \label{sec:decomposability}

We discuss a class of matrices for which the problem of determining whether two matrices lie in the same double coset reduces to similar problems involving smaller matrices.

\begin{prop}
  \label{theorem:decomposability}
  Suppose that $n=n_1+n_2$ and the the matrices $\alpha$ and $\alpha'$ have block matrix decompositions
  \begin{equation*}
    \alpha=\left(
    \begin{array}{cc}
      0 & \alpha_1\\
      \alpha_2 & 0
    \end{array}\right)
  \mbox{ and }
    \alpha'=\left(
    \begin{array}{cc}
      0 & \alpha_1'\\
      \alpha_2' & 0
    \end{array}\right),
  \end{equation*}
with $\alpha_i\in GL_{n_i}(A)$. Then $\alpha\sim \alpha'$ if and only if $\alpha_1\sim \alpha_1'$ and $\alpha_2\sim \alpha_2'$.
\end{prop}
\begin{proof}
  Write
  \begin{equation*}
        \left(
    \begin{array}{cc}
      b_1 & X\\
      0 & b_2
    \end{array}\right)
  \left(
    \begin{array}{cc}
      0 & \alpha_1\\
      \alpha_2 & 0
    \end{array}\right)
  \left(
    \begin{array}{cc}
      c_1 & Y\\
      0 & c_2
    \end{array}\right)
  =   \left(
    \begin{array}{cc}
      0 & \alpha_1'\\
      \alpha_2' & 0
    \end{array}\right)
  \end{equation*}
  where $b_1$ and $c_2$ are upper-triangular invertible $n_1\times n_1$-matrices, $b_2$ and $c_1$ are upper triangular invertible $n_2\times n_2$ matrices, $X$ is an $n_1\times n_2$ matrix and $Y$ is an $n_2\times n_1$ matrix.
  Multiplying out, and comparing the the lower-right blocks gives $b_2\alpha_2Y=0$, which, since $b_2$ and $\alpha_2$ are invertible, implies that $Y=0$.
Equating the remaining entries, and setting $Y=0$ gives that $b_2\alpha_2c_1=\alpha_2'$, $b_1\alpha_1c_2=\alpha_1'$ and $X=-b_1\alpha_1\alpha_2^{-1}$.
This shows that $\alpha\sim \alpha'$ if and only if $\alpha_1\sim \alpha_1'$ and $\alpha_2\sim \alpha_2'$.
\end{proof}

This proposition shows that the classification of the double cosets
for a given $n$ implies the classification for all $n'<n$. The
following corollary allows one to reduce the equivalence problem to
smaller $n$ for many matrices:
\begin{cor}
  \label{cor:decomposability}
   Suppose that $n=n_1+n_2$ and the the matrices $\alpha$ and $\alpha'$ have block matrix decompositions
  \begin{equation*}
    \alpha=\left(
    \begin{array}{cc}
      X_1 & \alpha_1\\
      \alpha_2 & X_2
    \end{array}\right)
  \mbox{ and }
    \alpha'=\left(
    \begin{array}{cc}
      X_1' & \alpha_1'\\
      \alpha_2' & X_2'
    \end{array}\right),
  \end{equation*}
with $\alpha_i\in GL_{n_i}(A)$. Then $\alpha\sim \alpha'$ if and only if $\alpha_1-X_1\alpha_2^{-1}X_2\sim \alpha_1'-X_2'\alpha_2^{\prime-1}X_2'$ and $\alpha_2\sim \alpha_2'$.
 \end{cor}
 \begin{proof}
   The upper triangular row and column operations described in Section~\ref{sec:upper_triangular} can be used to reduce
   $\alpha$ and $\alpha'$ to matrices of the type that occur in Proposition~\ref{theorem:decomposability},
   but with $\alpha_1$ and $\alpha_1'$ being replaced by $\alpha_1-X_1\alpha_2^{-1}X_2$ and $\alpha_1'-X_2'\alpha_2^{\prime-1}X_2'$ respectively.
 \end{proof}

\subsection{Intersection Invariants}
\label{sec:intersection}
Let $\F_0$ denote the \emph{standard flag} in $A^n$:
\begin{equation*}
  \F_0=(0=\F_0^0\subset \F_0^1\subset \cdots \subset \F_0^n=A^n)
\end{equation*}
where $\F_0^i$ is the $A$-module spanned by $\{\mathbf{e}_1,\ldots,\mathbf{e}_i\}$, and $\mathbf{e}_i$ is the $i$th standard basis vector in $A^n$.
$G$ acts transitively on the set of full free primitive flags over $A$ in $A^n$.
Thus the space of such flags is identified with $G/B$.
For $\alpha\in GL_n(A)$, consider the corresponding flag $\F=\alpha\F_0$ given by
\begin{equation*}
\F=(0=\F^0\subset\cdots\subset \F^n=A^n)=\alpha\F_0 \mbox{ where }\F^i=\alpha\F^i_0.
\end{equation*}
Clearly the isomorphism classes of the $A$-modules $\F^j\cap \F^i_0$
are invariants of the double coset. We will call these the {\em
intersection types}. The intersection types are related to the
column spaces of lower-left submatrices of $\alpha$ as follows: let
$[\alpha]^{ij}$ denote the lower-left $(n-i)\times j$ submatrix.
\begin{prop}
  \label{prop:matrix_to_type}
  The column space of $[\alpha]^{ij}$ is isomorphic to $\F^j/\F^j\cap \F^i_0$ as an $A$-module.
\end{prop}
\begin{proof}
  The map from $\F^j$ to the column space of $[\alpha]^{ij}$ is defined by taking the last $n-i$ entries of a vector. The kernel is clearly $\F^j\cap \F^i_0$.
\end{proof}

Furthermore, $\F$ induces a filtration on each
graded piece $\F^i_0/\F^{i-1}_0$ of the standard flag $\F_0$:
\begin{equation*}
  0=\frac{\F^0\cap \F^i_0}{\F^0\cap \F^{i-1}_0}\subset \frac{\F^1\cap \F^i_0}{\F^1\cap \F^{i-1}_0}\subset \cdots \subset \frac{\F^n\cap \F^i_0}{\F^n\cap \F^{i-1}_0}=\F^i_0/\F^{i-1}_0.
\end{equation*}
The $j$th graded piece of the above filtration is:
\begin{equation}
  \label{eq:type}
  \frac{\F^j \cap \F_0^i}{\F^{j-1}\cap \F_0^{i} + \F^j\cap \F_0^{i-1}}.
\end{equation}
Being a subquotient of $\F^i_0/\F^{i-1}_0\cong A$, it must be isomorphic to
$A/(p^{r_{ij}})$ for some $0\leq r_{ij}\leq k$ (if $k=\infty$, then
some of the $r_{ij}$'s will be infinite). The $B$-action on the
space $G/B$ of flags preserves the isomorphism classes of the
$A$-modules in \eqref{eq:type}.
Consider the matrix
$r(\alpha)=(r_{ij})$.
The above considerations show that it is
invariant under left and right multiplications in $B$, and that each
column sums to $k$.
A similar argument can be used to show that each
row sums to $k$.
We call $r(\alpha)$ \emph{the matrix of intersection numbers}
of $\alpha$.
When $k=1$, the matrix of intersection numbers is a permutation matrix, and is in fact, the unique permutation matrix that lies in the double coset of $\alpha$.
\emph{In this sense, the matrix of intersection numbers is a direct generalisation of the permutation associated to a matrix over a field by the Bruhat decomposition.}
\subsection{Permutation}
\label{sec:permutation} Consider the surjection induced by reduction
modulo $\wp$
\[
\mathbf{W}:B(A) \bsl \GL_n(A) / B(A) \longrightarrow B(\kk) \bsl
\GL_n(\kk) / B(\kk) \simeq S_n
\]
and view the double cosets as fibres over the field case. The image
of $\alpha\in \G$ in $B(\kk)\bsl \GL_n(\kk)/B(\kk)$ determines an
$n \times n$ permutation matrix $\mathbf{W}(\alpha)$ by the Bruhat
decomposition which is an invariant of the double coset of $\alpha$.
Given a permutation matrix $w$ of order $n$, let $N(w)$ denote the number of double cosets for which $\mathbf{W}=w$.
Say that the permutation $w$ is decomposable if there exists a partition $n=n_1+n_2$ and $w_1$ and $w_2$ permutation matrices of order $n_1$ and $n_2$ respectively such that
\begin{equation*}
  w=\left(
    \begin{array}{cc}
      0 & w_1\\
      w_2 & 0
    \end{array}\right)
\end{equation*}
as a block matrix.
By Corollary~\ref{cor:decomposability}, we have:
\begin{equation}
  \label{eq:decomposability}
  N(w)\geq N(w_2) \mbox{ for all } k.
\end{equation}

\section{The $n=2$ case}
\label{sec:n_is_2} In this case there are $k+1$ double cosets.
Geometrically, these double cosets parameterise the possible
intersections of two free primitive sub-modules of $A^2$ of rank
$1$. The intersection of two such submodules is a submodule of each
of them, and hence isomorphic to $A/\wp^r$ for some $0\leq r\leq k$.
We see that the intersection diagram, the intersection types and the
intersection numbers carry the same information; in fact they are
complete invariants. In terms of the permutation invariants, the
fibre over the trivial permutation consists of $k$ elements
corresponding to $r=1,\ldots,k$, and the fibre over the non-trivial
permutation consists of one element corresponding to $r=0$. In terms
of matrices, the set
\begin{equation*}
\left\{\left.\left(\begin{array}{cc} 1 & 0 \\
\pi^r & 1
\end{array}\right)~\right|~ 0 \le r \le k\right\}
\end{equation*}
 is a complete set of
representatives. The intersection numbers are given by
\begin{equation*}
    \alpha\in B
  \left(
    \begin{array}{cc}
      1 & 0 \\
      \pi^r & 1
    \end{array}
  \right)B
  \mbox{ if and only if }
      \mathbf{r}(\alpha)=  \left(
      \begin{array}{cc}
        r & k-r \\
        k-r & r
      \end{array}
    \right).
\end{equation*}

\section{The $n=3$ case}\label{sec:n_is_3}

\subsection{Fibration over the residue field} We have the following description of the fibres over permutations:
\begin{equation}\label{fib}
\begin{matrix}
B \bsl GL_3(A)/B \qquad &B \bsl M_2^{\bullet} / B & \val(\wp)\times
\val(\wp)& \val(\wp)\times \val(\wp)&
\val(\wp) & \val(\wp)& \{\star\} \\
 \downarrow&   \downarrow &
\downarrow & \downarrow & \downarrow & \downarrow &
\downarrow \\
B \bsl GL_3(\kk)/B \qquad &  \Bigl(\begin{smallmatrix} 1 & 0 & 0 \\
0 & 1 & 0 \\ 0 & 0 & 1 \end{smallmatrix}\Bigr)  &\Bigl(\begin{smallmatrix} 0 & 1 & 0 \\
1 & 0 & 0 \\ 0 & 0 & 1 \end{smallmatrix}\Bigr) & \Bigl(\begin{smallmatrix} 1 & 0 & 0 \\
0 & 0 & 1 \\ 0 & 1 & 0 \end{smallmatrix}\Bigr) & \Bigl(\begin{smallmatrix} 0 & 0 & 1 \\
1 & 0 & 0 \\ 0 & 1 & 0 \end{smallmatrix}\Bigr) & \Bigl(\begin{smallmatrix} 0 & 1 & 0 \\
0 & 0 & 1 \\ 1 & 0 & 0 \end{smallmatrix}\Bigr) & \Bigl(\begin{smallmatrix} 0 & 0 & 1 \\
0 & 1 & 0 \\ 1 & 0 & 0 \end{smallmatrix}\Bigr) \\
S_3  & 1 & s_1  &  s_2  &  s_1s_2 & s_2s_1 & s_1s_2s_1
\end{matrix}
\end{equation}
Except for the fibre over the trivial element (for which the notation used in the table is explained below), one easily verifies
that the fibres are indeed the ones written above. However, perhaps
for the case of the double cosets lying over the permutations
labelled $s_1$ and $s_2$ a remark is in order: any element
lying above $s_1$ can be brought to the form
\[
\left(\begin{matrix} 0 & 1 & 0 \\
1 & 0 & 0 \\ \pi^i & \pi^j & 1 \end{matrix}\right) ~, \qquad 1 \le
i,j \le k
\]
We see that these
lie in different double cosets by observing that they have different intersection
types. A similar argument holds for $s_2$.

As for the fibre over $\{1\}$, it is determined by the $2 \times 2$
lower left sub-matrix. We are therefore led to analyse the double
coset space of $2\times 2$ matrices $B \bsl M_2^{\bullet} /B$. Here
$M_2^{\bullet}$ denotes the set of those $2\times 2$ matrices for
which only the top right entry is a unit and (with a slight abuse of
notation) $B$ denotes the group of upper triangular matrices in
$GL_2(A)$.

\subsection{The space $B\bsl \M^\bullet / B$}
\label{sec:lowerleft}

For $\alpha \in \M^\bullet$, let $\val(\alpha)$ denote the valuation matrix.
The matrix is in {\em standard form} if the valuations of the non-zero entries form a standard tableaux (decreasing with respect to column numbers and increasing with respect to row numbers).
Every matrix in $M_2^\bullet$ can be reduced to standard form.
Assume that $\alpha$ is in standard form. Write
\begin{equation*}
  v(\alpha)=
  \left(
    \begin{array}{cc}
      i & 0 \\
      j & l
    \end{array}
  \right).
\end{equation*}
\begin{description}
\item [Discrete part] If at least one of the entries of $i$, $j$ or $l$ is $k$, the double coset is completely
determined by the valuation matrix, and the classification of the
orbits is given completely in terms of the matrix of intersection numbers.

\item [Non-discrete part] If none of the entries of $\val(\alpha)$ is equal to $k$, the matrix could be brought to the
form:
\[
\alpha(a)=\biggr(\begin{matrix} \pi^i & 1 \\ \pi^ja & \pi^l
\end{matrix}\biggl) 
\]
with $k>j>\max\{i,l\}\geq\min\{i,l\}> 0$ and $a \in A^{\times}$ (in fact can be taken in $A_{k-j}^{\times}$).\\
\end{description}
We shall consider the cases where $k$ is finite and infinite
simultaneously. When $k$ is finite, our equations are modulo a power
of $\pi$. When $k$ is infinite, then $A$ is a domain, and the
process of going modulo powers of $\pi$ should be ignored.
Let $\epsilon = \min \{j-i,j-l,i,l\}$ and $\delta(a)=\val(a-1)$.
\begin{prop}\label{BMB} Let $i,j$ and $l$ satisfy $k > j > \max \{i,l\} \geq \min \{i,l\} > 0$.
For every $a, a' \in A_{k-j}^{\times}$, let $\alpha(a)$ and $0 \leq
\delta(a) \leq k-j$ be as above. We have
\begin{enumerate}
\item If $j \neq i+l $, then
\[
\alpha(a) \sim \alpha(a') \Longleftrightarrow a \equiv a'\mod
\pi^{\min \{\epsilon, k-j\}}
\]
\item If $j = i+l$, then
\[
\alpha(a) \sim \alpha(a') \Longleftrightarrow \left\{
\begin{array}{ll}
    \delta(a)=\delta(a')=:\delta  \\
    a \equiv a'\mod \pi^{\min \{\epsilon+\delta, k-j\}} \\
\end{array}
\right.
\]
\end{enumerate}
\end{prop}

\begin{proof}
$\alpha(a) \sim \alpha(a')$ if and only if there exist $x_{11},x_{22},y_{11},y_{22}\in A^\times$ and $x_{12},y_{12}\in A$ such that
\begin{equation}\label{mat22}
\biggl(\begin{matrix} x_{11} & x_{12} \\ 0 & x_{22}
\end{matrix}\biggr) \biggl(\begin{matrix} \pi^i & 1 \\ \pi^ja & \pi^l
\end{matrix}\biggr) = \biggl(\begin{matrix} \pi^i & 1 \\ \pi^ja' & \pi^l
\end{matrix}\biggr) \biggl(\begin{matrix} y_{11} & y_{12} \\ 0 & y_{22}
\end{matrix}\biggr)
\end{equation}
Equating the entries on both sides of \eqref{mat22}
we get the following system of equations:
\begin{eqnarray*}
\pi^ix_{11}+\pi^jax_{12}& = & \pi^iy_{11}\\
x_{11}+\pi^lx_{12}&=& \pi^iy_{12}+y_{22}\\
\pi^jax_{22}&=&\pi^ja'y_{11}\\
\pi^lx_{22}&=& \pi^ja'y_{12}+\pi^ly_{22}
\end{eqnarray*}
and after successive substitutions of these equations we obtain
\begin{equation}\label{reduction}
(1-a'a\inv)y_{11}=(\pi^{j-i}a+\pi^{l})x_{12}+(\pi^{i}-\pi^{j-l}a')y_{12}
\mod \pi^{k-j}.
\end{equation}
The equation \eqref{reduction} can be solved if and only if the coefficient of $y_{11}$ on the left hand side lies in the ideal generated by the coefficients of $x_{12}$ and $y_{12}$ on the right hand side, in other words, if and only if
\begin{equation}\label{valuation}
v(a-a')\geq \min\{v(\pi^{j-i}a+\pi^{l}),v(\pi^{i}-a'\pi^{j-l}),
k-j\}.
\end{equation}
If $j \neq i+l$, then the right hand side of \eqref{valuation} is
$\min\{\epsilon,k-j\}$. This proves the first part of the
proposition. To prove the second part, we observe that if $j=i+l$,
then \eqref{valuation} becomes
\begin{equation}
\label{eq:valuation2} v(a-a')\geq \min\{l+v(a-1),i+v(a'-1),k-j\}.
\end{equation}
However, equating the valuations of the determinants of both sides
of \eqref{mat22} gives
\begin{equation*}
\min\{k-j,v(a-1)\}=\min\{k-j,v(a'-1)\},
\end{equation*}
and taking this condition into account allows us to replace $a'$
with $a$ in the right hand side of \eqref{valuation} to get
$\min\{k-j,i+v(a-1),l+v(a-1)\}$, which, when $j=i+l$, is the same as
$\min\{k-j,\epsilon+\delta\}$.

Conversely, it is easy to see that a solution of \eqref{reduction}
can always be extended to a solution of \eqref{mat22}.
\end{proof}

\subsection{Parametrisation of the double coset space}
\begin{cor}\label{parametrization}
The double coset space $B\bsl GL_3(A)/B$ is parameterised by
\[
B\bsl GL_3(A)/B = (B \bsl \M^{\bullet} /B) \coprod (\val(\wp)\times
\val(\wp)) \coprod (\val(\wp) \times \val(\wp)) \coprod \val(\wp)
\coprod \val(\wp) \coprod \{\star\}
\]
where
\[
\begin{split}
B \bsl \M^{\bullet} /B  = \coprod_{j=2}^{k-1} \left[
\left(\coprod_{\substack{1 \le i,l \le j-1 \\ i+l \ne j}}
A^{\times}_{\min\{\epsilon,k-j\}}\right) \coprod \left(
\coprod_{\substack{1 \le i,l \le j-1 \\ i+l = j}}\;
\coprod_{\delta=0}^{k-j}
A^{\times,\delta}_{\min\{\epsilon+\delta,k-j\}} \right) \right]
\\ \coprod \left\{\left.\left(\begin{matrix} i & 0 \\ j & l \end{matrix}\right)\right| [j=k, 1 \le i,l \le
k] ~\text{or}~ [i<j<k=l] ~\text{or}~ [l<j<k=i]\right\}
\end{split}
\]
and $A_i^{\times,\delta}=\{ a \in A_i^{\times} ~|~ v(a-1)=
\delta\}$.
\end{cor}

\begin{cor} If the residue field of $A$ is finite with $q$ elements then
\begin{equation*}
p_1(k) q^{\lfloor k/3 \rfloor} < |B\bsl GL_3(A_k)/B| < p_2(k)
q^{\lceil k/3 \rceil}
\end{equation*}
for some positive polynomials $p_1$ and $p_2$.
\end{cor}
\begin{proof}
For each permutation except for the trivial one, the number of
double cosets that lie over that permutation has polynomial growth
in $k$, hence does not affect the bounds. This is still the case for
the identity permutation, when at least one of $i$, $l$ or $j$ is
$k$, as there is only one double coset of the standard form of
\S\ref{sec:lowerleft} with valuations $i$, $l$ and $j$.

Otherwise, different powers of $q$ appear as the number of double
cosets with fixed values of $i$, $l$ and $j$. Here, there are two
possible cases:
\begin{align*}
&j \ne i+l  &\text{\# cosets}& = |A^{\times}_{\min\{\epsilon,k-j\}}|=q^{\min\{\epsilon,k-j\}}(1-q^{-1}) & &\quad \\
&j = i+l &\text{\# cosets}
&=|\coprod_{\delta=0}^{k-j}A_{\epsilon+\delta,
k-j}^{\times,\delta}|=\sum_{\delta=0}^{k-j}q^{\min\{\epsilon,
k-j-\delta\}} & &\quad
\end{align*}
Clearly, the highest power of $q$ that can occur in this way is
$\lceil k/3 \rceil$, hence the upper bound. This value is indeed
achieved when $i$, $l$ and $j/2$ are close to $k/3$, hence the lower
bound.
\end{proof}

\section{The $n=4$, $k=2$ case}\label{(4,2)}
We saw that when $n=3$ the number of double cosets does not depend
on the ring (e.g.,\ the characteristic of the residue field) when
$k=2$. This is no longer true when $n=4$. We illustrate this by
describing some of the double cosets lying over the permutation
\begin{equation*}
  w=
  \left(
    \begin{array}{cccc}
      0 & 1 & 0 & 0 \\
      1 & 0 & 0 & 0 \\
      0 & 0 & 0 & 1 \\
      0 & 0 & 1 & 0
    \end{array}
  \right)
\end{equation*}
Let $\sigma = \left(
  \begin{array}{cc}
    0 & 1\\
    1 & 0
  \end{array} \right)$, and $\tau(a)=\left(
  \begin{array}{cc}
    \pi & a\pi \\
    \pi & \pi
  \end{array}\right)$.
\begin{prop}
  \label{prop:42}
  For $a,a'\in A$, the (block) matrices
  \begin{equation*}
    \left(
      \begin{array}{cc}
        \sigma & 0 \\
        \tau(a) & \sigma
      \end{array}
    \right)
    \mbox{ and }
    \left(
      \begin{array}{cc}
        \sigma & 0 \\
        \tau(a') & \sigma
      \end{array}
    \right)
  \end{equation*}
lie in the same double coset if and only if $a\equiv a' \mbox{ mod } \pi$.
\end{prop}
\begin{proof}
  Suppose there exist invertible upper triangular $2\times 2$ matrices $D_1$, $D_2$, $C_1$ and $C_2$, and $2\times 2$ matrices $X$ and $Y$ such that
  \begin{equation*}
    \left(
      \begin{array}{cc}
        D_1 & X\\
        0 & D_2
      \end{array}
    \right)
    \left(
      \begin{array}{cc}
        \sigma & 0 \\
        \tau(a) & \sigma
      \end{array}
    \right)
    \left(
      \begin{array}{cc}
        C_1 & Y\\
        0 & C_2
      \end{array}
    \right)
    =
    \left(
      \begin{array}{cc}
        \sigma & 0 \\
        \tau(a') & \sigma
      \end{array}
    \right),
  \end{equation*}
then the following identities must hold:
\begin{eqnarray}
  \label{eq:topleft}
  D_1\sigma C_1+X \tau(a) C_1 & =& \sigma,\\
  \label{eq:topright}
  D_1\sigma Y + X \tau(a) Y + X\sigma C_2 & = & 0,\\
  \label{eq:bottomleft}
  D_2\tau(a)C_1 & = & \tau(a'),\\
  \label{eq:bottomright}
  D_2 \tau(a) Y + D_2\sigma C_2 & = &\sigma.
\end{eqnarray}
Now, reducing \eqref{eq:topleft} modulo $\pi$ gives
\begin{equation*}
  D_1\sigma \equiv \sigma C_1\inv \mod \pi.
\end{equation*}
Comparing the top-left and bottom-right entries of both sides shows
that $D_1$ and $C_1$ are in fact congruent to diagonal matrices
modulo $\pi$. Similarly, \eqref{eq:bottomright} can be used to show
that $D_2$ and $C_2$ are also congruent to diagonal matrices modulo
$\pi$. Thus $D_2$ can be written in the form $\big(
  \begin{smallmatrix}
    d_{11} & \pi d_{12}\\
    0 & d_{22}
  \end{smallmatrix} \big)$ where $d_{11}$ and $d_{22}$ are units.
Similarly, $C_1$ can be written in the form $\big(
  \begin{smallmatrix}
    c_{11} & \pi c_{12}\\
    0 & c_{22}
  \end{smallmatrix} \big)$ where $c_{11}$ and $c_{22}$ are units.
Substituting in \eqref{eq:bottomleft} and comparing entries gives
\begin{equation*}
  d_{11}c_{11}\equiv d_{22}c_{22}\equiv d_{22}c_{11}\equiv 1 \mbox{ mod }\pi,
\end{equation*}
\begin{equation*}
 d_{11}c_{22}a \equiv a' \mbox{ mod } \pi.
\end{equation*}
The above equations are readily seen to imply that $a\equiv a' \mbox{ mod } \pi$.
\end{proof}

\begin{cor}
  \label{theorem:42}
  When $n=4$ and $k=2$, there exist at least $|\kk|$ double cosets in $B\bsl G /B$ that lie above the permutation $w$.
\end{cor}
\noindent {\bf Remark.} Note that the intersection invariants do not
distinguish between all the double cosets. Among the matrices
considered in Proposition~\ref{prop:42}, the double coset
corresponding to $a\equiv 1 \mbox{ mod } \pi$ can be distinguished
from the others by the $(2,3)$th entry of the matrix of intersection
numbers described in \S\ref{sec:intersection}. However, the
intersection types and intersection numbers do not distinguish
between the remaining $(|\kk|-1)$ double cosets.
\section{The general picture}
\begin{theorem}
  \label{theorem:cases}
    The number of double cosets in $B\bsl GL_n(A)/B$ does not depend on $\kk$ if and only if $n\leq 2$ or $(n,k)=(3,2)$ or
  $k=1$.
\end{theorem}
\begin{proof}
  We have seen in Sections \ref{sec:n_is_2} and \ref{sec:n_is_3} that the number of double cosets does not depend on $\kk$ when $n\leq 2$ or when $n=3$ and $k=2$.
  We have also seen in Sections \ref{sec:n_is_3} and \ref{(4,2)} that for $n=3$ and $k>2$ and for $n=4$ and $k>1$ the number of double cosets does depend on $\kk$.
  In any other case $(n,k)$ there exists $n'<n$ such that the number of cosets depends on $\kk$ for the case $(n',k)$.
  Take a permutation matrix $w'$ of order $n'$ for which the number of double cosets with permutation invariant $w'$ in $GL_{n'}(A)$ depends on $\kk$.
  Let $w$ be the block matrix
  \begin{equation*}
    \left(
      \begin{array}{cc}
        0 & I_{n-n'} \\
        w' & 0
      \end{array}\right),
  \end{equation*}
  where $I_{n-n'}$ is the identity matrix of order $n-n'$.
  By Section~\ref{sec:permutation}, $N(w)\geq N(w')$ so that the number of double cosets with permutation invariant $w$ also depends on $\kk$.
\end{proof}
The complexity of the double coset space is sketched in the
following table
\[
\begin{array}{ccccccccccc}
    & \vline & 1 & 2 & 3 & 4 &   \cdots & n \\
   \hline
   1 & \vline & D & D & D & D &  \cdots  &  \\
   2 & \vline & D & D & D & N &  \cdots  &  \\
   3 & \vline & D & D & N & N & \cdots &    \\
   4 & \vline & D & D & N & N& \cdots &    \\
   5 & \vline & D & D & N & N& \cdots &   \\
   \vdots & \vline & \vdots & \vdots & \vdots & \vdots  &  &    \\
   k & \vline &  &  &  &  &  &
\end{array}
\]
Here $N$ stands for non-discrete and $D$ stands for discrete, to
indicate whether the double coset space depends or does not depend
(respectively) on the ring in question, e.g.,\ on the cardinality of
the residue field.

\providecommand{\bysame}{\leavevmode\hbox
to3em{\hrulefill}\thinspace}
\providecommand{\MR}{\relax\ifhmode\unskip\space\fi MR }
\providecommand{\MRhref}[2]{%
  \href{http://www.ams.org/mathscinet-getitem?mr=#1}{#2}
} \providecommand{\href}[2]{#2}

\vspace{\bigskipamount}
\vspace{\bigskipamount}
\vspace{\bigskipamount}

\begin{footnotesize}
\begin{quote}
Uri Onn\\
Einstein Institute of Mathematics,
Edmond Safra Campus, Givat Ram, \\
 Jerusalem 91904, Israel \\
{\tt urion@math.huji.ac.il} \\ \\

Amritanshu Prasad\\
The Institute of Mathematical Sciences, CIT campus, \\
Chennai 600 113, India \\
{\tt amri@imsc.res.in} \\ \\

Leonid Vaserstein\\
Department of Mathematics, Penn State University, University Park PA
\\ 16802-6401, USA \\
{\tt vstein@math.psu.edu} \\ \\
\end{quote}
\end{footnotesize}
\end{document}